%
%
%
%
%
%
\magnification=\magstephalf      
%
%
\vsize=7.5truein                 
\hsize=5.2truein                 
\newskip\stdskip                 
\stdskip=6pt plus3pt minus3pt    
\medskipamount=\stdskip          
\parindent=0pt                   
\parskip=\stdskip                
\abovedisplayskip=\stdskip       
\belowdisplayskip=\stdskip       
\mathsurround=0.75pt             
\overfullrule=0pt                
%
%
\def\ppar{\par\goodbreak\vskip 8pt plus 4pt minus 4pt}     
%
%
\def\stdspace{\hskip 0.75em plus 0.15em\ignorespaces}
\let\qua\stdspace 
%
%
%
%
%
%
%
\def\hexnumber#1{\ifcase#1 0\or 1\or 2\or 3\or 4\or 5\or 6\or 7\or 8\or
 9\or A\or B\or C\or D\or E\or F\fi}
%
%
\font\thirtnmsa=msam10 scaled 1315    
\font\tenmsa=msam10          \font\ninemsa=msam9
\font\sevenmsa=msam7         \font\sixmsa=msam6
\font\fivemsa=msam5
%
%
\newfam\msafam                  \textfont\msafam=\tenmsa
\scriptfont\msafam=\sevenmsa    \scriptscriptfont\msafam=\fivemsa
\edef\hexa{\hexnumber\msafam}        
\def\msa{\fam\msafam\tenmsa}         
%
%
\font\thirtnmsb=msbm10 scaled 1315   
\font\tenmsb=msbm10      \font\ninemsb=msbm9
\font\sevenmsb=msbm7     \font\sixmsb=msbm6
\font\fivemsb=msbm5
%
\newfam\msbfam                   \textfont\msbfam=\tenmsb       
\scriptfont\msbfam=\sevenmsb     \scriptscriptfont\msbfam=\fivemsb
\edef\hexb{\hexnumber\msbfam}    
\def\msb{\fam\msbfam\tenmsb}     
%
%
\font\thirtneufm=eufm10 scaled 1315   
\font\teneufm=eufm10                 \font\nineeufm=eufm9
\font\seveneufm=eufm7                \font\sixeufm=eufm6
\font\fiveeufm=eufm5
%
\newfam\eufmfam                    \textfont\eufmfam=\teneufm
\scriptfont\eufmfam=\seveneufm     \scriptscriptfont\eufmfam=\fiveeufm
\edef\hexf{\hexnumber\eufmfam}      
\def\frak{\fam\eufmfam\teneufm}     
%
%
%
\font\thirtnrm=cmr10 scaled 1315    
\font\ninerm=cmr9                   \font\sixrm=cmr6   
%
\font\thirtni=cmmi10 scaled 1315    
\font\ninei=cmmi9                   \font\sixi=cmmi6  
%
\font\thirtnsy=cmsy10 scaled 1315   
\font\ninesy=cmsy9                  \font\sixsy=cmsy6  
%
\font\thirtnbf=cmbx10 scaled 1315   
\font\ninebf=cmbx9                  \font\sixbf=cmbx6  
%
%
\font\thirtnex=cmex10 scaled 1315   
\font\nineex=cmex9                  
%
%
\font\thirtnit=cmti10 scaled 1315  
\font\nineit=cmti9                  
%
\font\thirtnsl=cmsl10 scaled 1315  
\font\ninesl=cmsl9                  
%
\font\thirtntt=cmtt10 scaled 1315  
\font\ninett=cmtt9                  
%
%
%
%
\def\small{%
%
%
\textfont0=\ninerm \scriptfont0=\sixrm \scriptscriptfont0=\fiverm
\def\rm{\fam0\ninerm}
%
%
\textfont1=\ninei \scriptfont1=\sixi \scriptscriptfont1=\fivei
%
%
\textfont2=\ninesy \scriptfont2=\sixsy \scriptscriptfont2=\fivesy
%
%
\textfont3=\nineex \scriptfont3=\nineex \scriptscriptfont3=\nineex
%
%
\textfont\bffam=\ninebf \scriptfont\bffam=\sixbf
\scriptscriptfont\bffam=\fivebf \def\bf{\fam\bffam\ninebf}%
%
%
\textfont\itfam=\nineit \def\it{\fam\itfam\nineit}%
\textfont\slfam=\ninesl \def\sl{\fam\slfam\ninesl}%
\textfont\ttfam=\ninett \def\tt{\fam\ttfam\ninett}%
%
%
%
\textfont\msafam=\ninemsa \scriptfont\msafam=\sixmsa
\scriptscriptfont\msafam=\fivemsa \def\msa{\fam\msafam\ninemsa}%
%
%
\textfont\msbfam=\ninemsb \scriptfont\msbfam=\sixmsb
\scriptscriptfont\msbfam=\fivemsb \def\msb{\fam\msbfam\ninemsb}%
%
%
\textfont\eufmfam=\nineeufm  \scriptfont\eufmfam=\sixeufm
\scriptscriptfont\eufmfam=\fiveeufm \def\frak{\fam\eufmfam\nineeufm}%
%
%
%
\normalbaselineskip=11pt%
\setbox\strutbox=\hbox{\vrule height8pt depth3pt width0pt}%
%
%
\normalbaselines\rm
%
%
\stdskip=4pt plus2pt minus2pt    
\medskipamount=\stdskip          
\parskip=\stdskip                
\abovedisplayskip=\stdskip       
\belowdisplayskip=\stdskip       
\def\ppar{\par\goodbreak\vskip 6pt plus 3pt minus 3pt}%
%
%
\def\section##1{\global\advance\sectionnumber by 1
\vskip-\lastskip\penalty-800\vskip 20pt plus10pt minus5pt 
\egroup{\bf\number\sectionnumber\quad##1}\bgroup\small         
\vskip 6pt plus3pt minus3pt
\nobreak\resultnumber=1}
}    
%
\def\beginsmall{\bgroup\small}
\let\endsmall\egroup
%
%
%
%
\def\large{%
\textfont0=\thirtnrm \scriptfont0=\ninerm \scriptscriptfont0=\sevenrm
\def\rm{\fam0\thirtnrm}%
\textfont1=\thirtni \scriptfont1=\ninei \scriptscriptfont1=\seveni
\textfont2=\thirtnsy \scriptfont2=\ninesy \scriptscriptfont2=\sevensy
\textfont3=\thirtnex \scriptfont3=\thirtnex \scriptscriptfont3=\thirtnex
\textfont\bffam=\thirtnbf \scriptfont\bffam=\ninebf
\scriptscriptfont\bffam=\sevenbf \def\bf{\fam\bffam\thirtnbf}%
\textfont\itfam=\thirtnit \def\it{\fam\itfam\thirtnit}%
\textfont\slfam=\thirtnsl \def\sl{\fam\slfam\thirtnsl}%
\textfont\ttfam=\thirtntt \def\tt{\fam\ttfam\thirtntt}%
\textfont\msafam=\thirtnmsa \scriptfont\msafam=\ninemsa
\scriptscriptfont\msafam=\sevenmsa \def\msa{\fam\msafam\thirtnmsa}%
\textfont\msbfam=\thirtnmsb \scriptfont\msbfam=\ninemsb
\scriptscriptfont\msbfam=\sevenmsb \def\msb{\fam\msbfam\thirtnmsb}%
\textfont\eufmfam=\thirtneufm  \scriptfont\eufmfam=\nineeufm
\scriptscriptfont\eufmfam=\seveneufm \def\frak{\fam\eufmfam\teneufm}%
\normalbaselineskip=16pt%
\setbox\strutbox=\hbox{\vrule height11.5pt depth4.5pt width0pt}%
\normalbaselines\rm}%
\let\Large\large   
%

%

%
\mathchardef\plussquare="0\hexa01
\mathchardef\nge="3\hexb0B
\mathchardef\maltesecross="0\hexa7A
\mathchardef\del="0\hexf01
%
%
%
%
\font\sc=cmcsc10
%
%
%
%
\def\sqr#1#2{{\vcenter{\vbox{\hrule  height.#2truept
	\hbox{\vrule width.#2truept height#1truept 
	\kern#1truept \vrule width.#2truept}
	\hrule height.#2truept}}}}
\def\sq{\sqr55}    
%
%
%
%
\newcount\sectionnumber            
\newcount\resultnumber             
\sectionnumber=0\resultnumber=1    
%
%
%
\def\section#1{\global\advance\sectionnumber by 1
\xdef\nextkey{\number\sectionnumber}
\vskip-\lastskip\penalty-800\vskip 20pt plus10pt minus5pt 
{\large\bf\number\sectionnumber\quad#1}         
\vskip 8pt plus4pt minus4pt
\nobreak\resultnumber=1}      
%
%
%
%
%
\def\sh#1{\vskip-\lastskip\ppar{\bf #1}\par\nobreak\medskip}         
%
%
%
%

%
\def\proc#1{\xdef\nextkey{\number\sectionnumber.\number\resultnumber}%
\vskip-\lastskip\ppar\bf%
\noindent#1\ \number\sectionnumber.\number\resultnumber
\stdspace\sl\global\advance\resultnumber by 1\ignorespaces}
\def\endproc{\rm\ppar} 
%
%
\def\prf{\vskip-\lastskip\ppar\noindent{\bf Proof}%
\stdspace\rm}                            
\def\endprf{\unskip\stdspace\hbox{}
\hfill$\sq$\par\medskip}                 
\def\proof#1{\vskip-\lastskip\ppar\noindent{\bf#1}%
\stdspace\rm\ignorespaces}        
%
%
%
%
%
%
%
%
\def\proclaim#1{\vskip-\lastskip\ppar\bf%
\noindent#1\stdspace\sl\ignorespaces} 
\let\endproclaim\endproc
%
%
%
%
\def\rk#1{\vskip-\lastskip\ppar{\bf #1}\stdspace\ignorespaces}                

%
%
%
%
%
%
\def\label{\xdef\nextkey{\number\sectionnumber.\number\resultnumber}%
\number\sectionnumber.\number\resultnumber
\global\advance\resultnumber by 1}
%
%
%
%
%
%
%
%
%
%
%
%
%
%
%
%
\newcount\refnumber              
\refnumber=1                     
\long\def\reflist#1\endreflist{%
\long\def\thereflist{#1}{\def\refkey##1##2\par{\xdef##1{\number\refnumber}%
\global\advance\refnumber by 1}%
\def\key##1##2\par{\expandafter\xdef%
\csname##1\endcsname{\number\refnumber}%
\global\advance\refnumber by 1}#1\par}}
\long\def\references{%
\penalty-800\vskip-\lastskip\vskip 15pt plus10pt minus5pt 
{\large\bf References}\ppar 
{\leftskip=25pt\frenchspacing    
\small\parskip=3pt plus2pt       
\def\refkey##1##2\par{\noindent  
\llap{[##1]\stdspace}\ignorespaces##2\par}         
\def\key##1##2\par{\noindent  
\llap{[\ref{##1}]\stdspace}\ignorespaces##2\par}  
\def\,{\thinspace}\thereflist\par}}
%
%
%
\newcount\footnotenumber         
\footnotenumber=1                
\def\fnote#1{\xdef\nextkey{\number\footnotenumber}%
{\small\ifnum\footnotenumber>9\parindent=14pt%
\else\parindent=10pt\fi\footnote{$^{\number\footnotenumber}$}%
{\hglue-5pt#1}\global\advance\footnotenumber by 1}}
%
%
%
%
%
%
%
\newcount\figurenumber          
\figurenumber=1                 
\def\caption#1{\xdef\nextkey{\number\figurenumber}%
\cl{\small Figure \number\figurenumber: #1}%
\global\advance\figurenumber by 1}
\def\figurelabel{\xdef\nextkey{\number\figurenumber}%
\cl{\small Figure \number\figurenumber}%
\global\advance\figurenumber by 1}
\long\def\figure#1\endfigure{{\xdef\nextkey{\number\figurenumber}%
\let\captiontext\relax\def\caption##1{\xdef\captiontext{##1}}%
\midinsert\cl{\ignorespaces#1\unskip\unskip\unskip\unskip}\vglue6pt\cl{\small 
Figure \number\figurenumber\ifx\captiontext\relax\else: \captiontext
\fi}\endinsert\global\advance\figurenumber by 1}}
%
%
%
%
%
%
%
\def\nextkey{??}   
%
\def\key#1{\expandafter\xdef\csname #1\endcsname{\nextkey}}
\def\ref#1{\expandafter\ifx\csname #1\endcsname\relax
\immediate\write16{Reference {#1} undefined}??\else
\csname #1\endcsname\fi}
%
%
%
%
%
%
%
\newread\gtinfile
\newwrite\gtreffile
\def\useforwardrefs{
\openin\gtinfile\jobname.ref
\ifeof\gtinfile
\closein\gtinfile
\immediate\write16{No file \jobname.ref}
\else
\closein\gtinfile
\input \jobname.ref
\fi
\immediate\openout\gtreffile \jobname.ref
%
%
\def\key##1{{\def\\{\noexpand}%
\expandafter\xdef\csname ##1\endcsname{\nextkey}%
\immediate\write\gtreffile{\\\expandafter\\\def\\\csname ##1\\\endcsname%
{\nextkey}}}}
%
%
\long\def\reflist##1\endreflist{%
\long\def\thereflist{##1}{\def\refkey####1####2\par{\xdef####1{%
\number\refnumber}{\def\\{\noexpand}\immediate\write\gtreffile
{\\\def\\####1{\number\refnumber}}}\global\advance\refnumber by 1}%
\def\key####1####2\par{\expandafter\xdef%
\csname####1\endcsname{\number\refnumber}%
{\def\\{\noexpand}\immediate\write\gtreffile
{\\\expandafter\\\def\\\csname ####1\\\endcsname{\number\refnumber}}}
\global\advance\refnumber by 1}##1\par}}
\long\def\biblio##1\endbiblio{\reflist##1\endreflist\references}%
%
%
\def\numkey##1{{\def\\{\noexpand}%
\xdef##1{\number\sectionnumber.\number\resultnumber}
\immediate\write\gtreffile{\\\def\\##1%
{\number\sectionnumber.\number\resultnumber}}}}
\def\seckey##1{{\def\\{\noexpand}\xdef##1{\number\sectionnumber}
\immediate\write\gtreffile{\\\def\\##1{\number\sectionnumber}}}}
\def\figkey##1{\xdef##1{\number\figurenumber}%
{\def\\{\noexpand}\immediate\write\gtreffile%
{\\\def\\##1{\number\figurenumber}}}
\number\figurenumber\global\advance\figurenumber by 1}
}   
%
%
%
%
\def\figkey#1{\xdef#1{\number\figurenumber}%
\number\figurenumber\global\advance\figurenumber by 1}
\def\fig#1#2\endfig{%
\midinsert\cl{#2}\vglue6pt\cl{\small Figure #1}\endinsert}
\def\newfig{\number\figurenumber\global\advance\figurenumber by 1}
\def\numkey#1{\xdef#1{\number\sectionnumber.\number\resultnumber}}
\def\seckey#1{\xdef#1{\number\sectionnumber}}
%
%
%
%
%
%
%
%
%
\def\verb{\catcode`\"=\active}       
\def\brev{\catcode`\"=12}            
\brev                                
\verb                                
{\obeyspaces\gdef {\ }}              
{\catcode`\`=\active\gdef`{\relax\lq}}
\def"{%
\begingroup\baselineskip=12pt\def\par{\leavevmode\endgraf}%
\tt\obeylines\obeyspaces\parskip=0pt\parindent=0pt%
\catcode`\$=12\catcode`\&=12\catcode`\^=12\catcode`\#=12%
\catcode`\_=12\catcode`\~=12%
\catcode`\{=12\catcode`\}=12\catcode`\%=12\catcode`\\=12%
\catcode`\`=\active\let"\endgroup}
\brev      
%
%
%
%
%
%
\def\items{\par\leftskip = 25pt}           
\def\enditems{\par\leftskip = 0pt}         
\def\item#1{\par\leavevmode\llap{#1\stdspace}%
\ignorespaces}                             
%
%

%
%
\def\np{\vfil\eject}         
\def\nl{\hfil\break}         
\def\cl{\centerline}         
\def\agt{{\mathsurround=0pt\it$\cal A\mskip-.7mu$lgebraic \&\ 
$\cal G\mskip-2mu$eometric $\cal T\!\!$opology}}  
%
%
%

%
%
%
%
%
\def\title#1{\def\thetitle{#1}}

\def\author#1{\edef\previousauthors{\theauthors}
 \ifx\theauthors\relax\def\theauthors{#1}\else
 \def\theauthors{\previousauthors\par#1}\fi}

\let\authors\author        
\def\address#1{\edef\previousaddresses{\theaddress}
 \ifx\theaddress\relax\def\theaddress{#1}\else
 \def\theaddress{\previousaddresses\par\vskip 2pt\par#1}\fi}
\def\secondaddress#1{\edef\previousaddresses{\theaddress}
 \ifx\theaddress\relax\def\theaddress{#1}\else
 \def\theaddress{\previousaddresses\par{\rm and}\par#1}\fi}   

\def\email#1{\edef\previousemails{\theemail}
 \ifx\theemail\relax\def\theemail{#1}\else
 \def\theemail{\previousemails\hskip 0.75em\relax#1}\fi}
\def\secondemail#1{\edef\previousemails{\theemail}
 \ifx\theemail\relax\def\theemail{#1}\else
 \def\theemail{\previousemails\hskip 0.75em{\rm and}\hskip 0.75em
 \relax#1}\fi}
\def\url#1{\edef\previousurls{\theurl}
 \ifx\theurl\relax\def\theurl{#1}\else
 \def\theurl{\previousurls\hskip 0.75em\relax#1}\fi}
\def\secondurl#1{\edef\previousurls{\theurl}
 \ifx\theurl\relax\def\theurl{#1}\else
 \def\theurl{\previousurls\hskip 0.75em{\rm and}\hskip 0.75em
 \relax#1}\fi}
\long\def\abstract#1\endabstract{\long\def\theabstract{#1}}
\def\primaryclass#1{\def\theprimaryclass{#1}}
\def\secondaryclass#1{\def\thesecondaryclass{#1}}
\def\keywords#1{\def\thekeywords{#1}}
%
%
\let\\\par\let\thetitle\relax\let\theshorttitle\relax
\let\theauthors\relax\let\theshortauthors\relax
\let\theaddress\relax\let\theshortaddress\relax
\let\theemail\relax\let\theurl\relax
\let\theabstract\relax\let\theprimaryclass\relax
\let\thesecondaryclass\relax\let\thekeywords\relax
%
%
%
%
\long\def\maketitlepage{    

\vglue 0.2truein   

%
{\parskip=0pt\leftskip 0pt plus 1fil\def\\{\par\smallskip}{\large
\bf\thetitle}\par\medskip}   

\vglue 0.15truein 

%
{\parskip=0pt\leftskip 0pt plus 1fil\def\\{\par}{\sc\theauthors}
\par\medskip}%
 
\vglue 0.1truein 

%
{\small\parskip=0pt
{\leftskip 0pt plus 1fil\def\\{\par}{\sl\theaddress}\par}
\ifx\theemail\relax\else  
\vglue 5pt \def\\{\stdspace{\rm and}\stdspace} 
\cl{Email:\stdspace\tt\theemail}\fi
\ifx\theurl\relax\else    
\vglue 5pt \def\\{\stdspace{\rm and}\stdspace} 
\cl{URL:\stdspace\tt\theurl}\fi\par}

\vglue 7pt 

{\bf Abstract}

\vglue 5pt

\theabstract

\vglue 7pt 

{\bf AMS Classification numbers}\quad Primary:\quad \theprimaryclass\par

Secondary:\quad \thesecondaryclass

\vglue 5pt 

{\bf Keywords:}\quad \thekeywords

\np  

}    
%
%
\long\def\makeshorttitle{    


%
{\parskip=0pt\leftskip 0pt plus 1fil\def\\{\par\smallskip}{\large
\bf\thetitle}\par\medskip}   

\vglue 0.05truein 

%
{\parskip=0pt\leftskip 0pt plus 1fil\def\\{\par}{\sc\theauthors}
\par\medskip}%
 
\vglue 0.03truein 

%
{\small\parskip=0pt
{\leftskip 0pt plus 1fil\def\\{\par}{\sl\ifx\theshortaddress\relax
\theaddress\else\theshortaddress\fi}\par}
\ifx\theemail\relax\else  
\vglue 5pt \def\\{\stdspace{\rm and}\stdspace} 
\cl{Email:\stdspace\tt\theemail}\fi
\ifx\theurl\relax\else    
\vglue 5pt \def\\{\stdspace{\rm and}\stdspace} 
\cl{URL:\stdspace\tt\theurl}\fi\par}

\vglue 10pt 


{\small\leftskip 25pt\rightskip 25pt{\bf Abstract}\stdspace\theabstract

{\bf AMS Classification}\stdspace\theprimaryclass
\ifx\thesecondaryclass\relax\else; \thesecondaryclass\fi\par
{\bf Keywords}\stdspace \thekeywords\par}
\vglue 7pt
}    
\let\maketitle\makeshorttitle        
%
%

\def\volumenumber#1{\def\thevolumenumber{#1}}
\def\volumeyear#1{\def\thevolumeyear{#1}}
\def\pagenumbers#1#2{\def\startpage{#1}\def\finishpage{#2}}
\def\published#1{\def\publishdate{#1}}
\def\received#1{\def\receiveddate{#1}}
\def\revised#1{\def\reviseddate{#1}}
\let\reviseddate\relax
\volumenumber{X}
\volumeyear{20XX}
\pagenumbers{1}{XXX}
\published{XX Xxxember 20XX}

\long\def\makeagttitle{   
\agt\hfill      
\hbox to 60truept{\vbox to 0pt{\vglue -14truept{\bf [Logo here]}\vss}\hss}
\break
{\small Volume \thevolumenumber\ (\thevolumeyear)
\startpage--\finishpage\nl
Published: \publishdate}

\vglue .2truein

{\parskip=0pt\leftskip 0pt plus 1fil\def\\{\par\smallskip}{\large
\bf\thetitle}\par\medskip}   
\vglue 0.05truein 

%
{\parskip=0pt\leftskip 0pt plus 1fil\def\\{\par}{\sc\theauthors}
\par\medskip}%
 
\vglue 0.03truein 


{\small\leftskip 25truept\rightskip 25truept{\bf Abstract}\stdspace\theabstract

{\bf AMS Classification}\stdspace\theprimaryclass
\ifx\thesecondaryclass\relax\else; \thesecondaryclass\fi\par
{\bf Keywords}\stdspace \thekeywords\par}\vglue 7truept

}   


\def\Addresses{\bigskip
{\small \parskip 0pt \leftskip 0pt \rightskip 0pt plus 1fil \def\\{\par}
\sl\theaddress\par\medskip \rm Email:\stdspace\tt\theemail\par
\ifx\theurl\relax\else\smallskip \rm URL:\stdspace\tt\theurl\par\fi}}

\def\agtart{
\hoffset 14truemm
\voffset 31truemm
\font\phead=cmsl9 scaled 950
\font\pnum=cmbx10 scaled 913
\font\pfoot=cmsl9 scaled 950
\headline{\vbox to 0pt{\vskip -4.5mm\line{\small\phead\ifnum
\count0=\startpage ISSN numbers are printed here
\hfill {\pnum\folio}\else\ifodd\count0\def\\{ }%
\ifx\theshorttitle\relax\thetitle\else\theshorttitle\fi\hfill{\pnum\folio}
\else\def\\{ and }{\pnum\folio}\hfill\ifx\theshortauthors\relax\theauthors
\else\theshortauthors\fi\fi\fi}\vss}}
\footline{\vbox to 0pt{\vglue 0mm\line{\small\pfoot\ifnum\count0=\startpage
Copyright declaration is printed here\hfill\else
\agt, Volume \thevolumenumber\ (\thevolumeyear)\hfill\fi}\vss}}
\let\maketitle\makeagttitle\let\makeshorttitle\makeagttitle}


\def\ifplaintex{\expandafter\ifx\csname documentclass\endcsname\relax}

\def\gtp{{\mathsurround=0pt\it $\cal G\mskip-2mu$eometry \&\ 
$\cal T\!\!$opology $\cal P\!$ublications}}  

\def\recd{{\small Received:\qua\receiveddate\ifx\reviseddate\relax
\else\qquad Revised:\qua\reviseddate\fi\par}} 


\def\lognumber#1{\def\thelognumber{#1}}
\def\volumenumber#1{\def\thevolumenumber{#1}}
\def\volumeyear#1{\def\thevolumeyear{#1}}
\def\papernumber#1{\def\thepapernumber{#1}}
\def\pagenumbers#1#2{\def\startpage{#1}\def\finishpage{#2}}
\def\published#1{\def\publishdate{#1}}

\def\received#1{\def\receiveddate{#1}}
\def\revised#1{\def\reviseddate{#1}}
\def\accepted#1{\def\accepteddate{#1}}

\def\asciiauthors#1{\def\theasciiauthors{#1}}
\def\asciiaddress#1{\def\theasciiaddress{#1}}

\def\coverauthors#1{\def\thecoverauthors{#1}}


\let\\\par\let\thelognumber\relax\let\thevolumenumber\relax
\let\thepapernumber\relax\let\thevolumeyear\relax\let\startpage\relax
\let\finishpage\relax\let\publishdate\relax\let\receiveddate\relax
\let\reviseddate\relax\let\accepteddate\relax\let\theasciititle\relax
\let\theasciiauthors\relax\let\theasciiaddress\relax
\let\theasciiabstract\relax

\let\thecoverauthors\relax\let\theasciiemail\relax


\ifplaintex
\font\logobig=cmssbx10 scaled 3836
\font\logomed=cmssbx10 scaled 2557
\else
\font\logobig=cmssbx10 scaled 4200
\font\logomed=cmssbx10 scaled 2800
\fi

\long\def\makeagttitle{   
\count0=\startpage
\agt\hfill      
\hbox to 45truept{\vbox to 0pt{\vglue -13truept{\logomed A\kern -.37em{\logobig 
T}\kern -.38em G}\vss}\hss}
\break
{\small Volume \thevolumenumber\ (\thevolumeyear)
\startpage--\finishpage\nl
Published: \publishdate}

\vglue .25truein

{\parskip=0pt\leftskip 0pt plus
1fil\def\\{\par\smallskip}{\Large\bf\thetitle}\par\medskip} \vglue
0.05truein

%
{\parskip=0pt\leftskip 0pt plus 1fil\def\\{\par}{\sc\theauthors}
\par\medskip}%
 
\vglue 0.03truein 


{\small\leftskip 25truept\rightskip 25truept{\bf Abstract}\stdspace\theabstract

{\bf AMS Classification}\stdspace\theprimaryclass
\ifx\thesecondaryclass\relax\else; \thesecondaryclass\fi\par
{\bf Keywords}\stdspace \thekeywords\par}\vglue 7truept

}   

\ifplaintex
\hoffset 14truemm
\voffset 31truemm
\font\phead=cmsl9 scaled 950
\font\pnum=cmbx10 scaled 913
\font\pfoot=cmsl9 scaled 950
\headline{\vbox to 0pt{\vskip -4.5mm\line{\small\phead\ifnum
\count0=\startpage ISSN 1472-2739 (on-line) 1472-2747 (printed)
\hfill {\pnum\folio}\else\ifodd\count0\def\\{ }%
\ifx\theshorttitle\relax\thetitle\else\theshorttitle\fi\hfill{\pnum\folio}
\else\def\\{ and }{\pnum\folio}\hfill\ifx\theshortauthors\relax\theauthors
\else\theshortauthors\fi\fi\fi}\vss}}
\footline{\vbox to 0pt{\vglue 0mm\line{\small\pfoot\ifnum\count0=\startpage
\copyright\ \gtp\hfill\else
\agt, Volume \thevolumenumber\ (\thevolumeyear)\hfill\fi}\vss}}
\else
\headsep 23pt
\footskip 35pt
\hoffset -4truemm
\voffset 12.5truemm
\font\lhead=cmsl9 scaled 1050
\font\lnum=cmbx10 
\font\lfoot=cmsl9 scaled 1050
\makeatletter
\def\@oddhead{{\small\lhead\ifnum\count0=\startpage ISSN 1472-2739 
(on-line) 1472-2747 (printed)\hfill {\lnum\number\count0}\else\ifodd\count0
\def\\{ }\ifx\theshorttitle\relax \thetitle \else\theshorttitle\fi\hfill
{\lnum\number\count0}\else\def\\{ and }{\lnum\number\count0}
\hfill\ifx\theshortauthors\relax 
\theauthors\else\theshortauthors\fi\fi\fi}}\def\@evenhead{\@oddhead}
\def\@oddfoot{\small\lfoot\ifnum\count0=\startpage\copyright\ \gtp\hfill\else
\agt, Volume \thevolumenumber\ (\thevolumeyear)\hfill\fi}
\def\@evenfoot{\@oddfoot}
\makeatother
\fi
\let\maketitlepage\makeagttitle
\let\makeshorttitle\maketitlepage
\let\maketitle\maketitlepage


\newwrite\gtoutfile
\long\gdef\makeheadfile{  
{\def\\{, }\def\s{ }
\immediate\openout\gtoutfile head.xxx
\immediate\write\gtoutfile{To: math@arxiv.org}
\immediate\write\gtoutfile{Subject: put OR rep NNNNN:ppppp}
\immediate\write\gtoutfile{--text follows this line--}
\immediate\write\gtoutfile{Proxy-for: \ifx\theasciiauthors\relax
\theauthors\else\theasciiauthors\fi\s<\ifx\theasciiemail\relax\theemail\else\theasciiemail\fi>}
\immediate\write\gtoutfile{\noexpand\\}
\immediate\write\gtoutfile{Authors: \ifx\theasciiauthors\relax
\theauthors\else\theasciiauthors\fi}
{\def\\{ }\immediate\write\gtoutfile{Title: \ifx\theasciititle\relax
\thetitle\else\theasciititle\fi}}
\immediate\write\gtoutfile{Subj-class: GT or SG, GR etc}
\immediate\write\gtoutfile{MSC-class: \theprimaryclass\ifx\thesecondaryclass\relax\else, \thesecondaryclass\fi}
\immediate\write\gtoutfile{Journal-ref: Algebr. Geom. Topol. \thevolumenumber\s
(\thevolumeyear) \startpage-\finishpage}
\immediate\write\gtoutfile{Comments: Published by Algebraic and
Geometric Topology at}
\immediate\write\gtoutfile{\s\s\s  http://www.maths.warwick.ac.uk/agt/AGTVol\thevolumenumber/agt-\thevolumenumber-\thepapernumber.abs.html}
\immediate\write\gtoutfile{\noexpand\\}
\immediate\write\gtoutfile{}
\ifx\theasciiabstract\relax
\immediate\write\gtoutfile{\theabstract}\else
\immediate\write\gtoutfile{\theasciiabstract}\fi
\immediate\write\gtoutfile{}
\immediate\write\gtoutfile{\noexpand\\}
\immediate\write\gtoutfile{}
\immediate\closeout\gtoutfile}}  

\def\maketitlepage{\makeagttitle\makeheadfile}
\let\makeshorttitle\maketitlepage
\let\maketitle\maketitlepage


\def\ifplaintex{\expandafter\ifx\csname documentclass\endcsname\relax}

\def\gtp{{\mathsurround=0pt\it $\cal G\mskip-2mu$eometry \&\ 
$\cal T\!\!$opology $\cal P\!$ublications}}  

\def\recd{{\small Received:\qua\receiveddate\ifx\reviseddate\relax
\else\qquad Revised:\qua\reviseddate\fi\par}} 


\def\lognumber#1{\def\thelognumber{#1}}
\def\volumenumber#1{\def\thevolumenumber{#1}}
\def\volumeyear#1{\def\thevolumeyear{#1}}
\def\papernumber#1{\def\thepapernumber{#1}}
\def\pagenumbers#1#2{\def\startpage{#1}\def\finishpage{#2}}
\def\published#1{\def\publishdate{#1}}

\def\received#1{\def\receiveddate{#1}}
\def\revised#1{\def\reviseddate{#1}}
\def\accepted#1{\def\accepteddate{#1}}

\def\asciiauthors#1{\def\theasciiauthors{#1}}
\def\asciiaddress#1{\def\theasciiaddress{#1}}

\def\coverauthors#1{\def\thecoverauthors{#1}}


\let\\\par\let\thelognumber\relax\let\thevolumenumber\relax
\let\thepapernumber\relax\let\thevolumeyear\relax\let\startpage\relax
\let\finishpage\relax\let\publishdate\relax\let\receiveddate\relax
\let\reviseddate\relax\let\accepteddate\relax\let\theasciititle\relax
\let\theasciiauthors\relax\let\theasciiaddress\relax
\let\theasciiabstract\relax

\let\thecoverauthors\relax\let\theasciiemail\relax


\ifplaintex
\font\logobig=cmssbx10 scaled 3836
\font\logomed=cmssbx10 scaled 2557
\else
\font\logobig=cmssbx10 scaled 4200
\font\logomed=cmssbx10 scaled 2800
\fi

\long\def\makeagttitle{   
\count0=\startpage
\agt\hfill      
\hbox to 45truept{\vbox to 0pt{\vglue -13truept{\logomed A\kern -.37em{\logobig 
T}\kern -.38em G}\vss}\hss}
\break
{\small Volume \thevolumenumber\ (\thevolumeyear)
\startpage--\finishpage\nl
Published: \publishdate}

\vglue .25truein

{\parskip=0pt\leftskip 0pt plus
1fil\def\\{\par\smallskip}{\Large\bf\thetitle}\par\medskip} \vglue
0.05truein

%
{\parskip=0pt\leftskip 0pt plus 1fil\def\\{\par}{\sc\theauthors}
\par\medskip}%
 
\vglue 0.03truein 


{\small\leftskip 25truept\rightskip 25truept{\bf Abstract}\stdspace\theabstract

{\bf AMS Classification}\stdspace\theprimaryclass
\ifx\thesecondaryclass\relax\else; \thesecondaryclass\fi\par
{\bf Keywords}\stdspace \thekeywords\par}\vglue 7truept

}   

\ifplaintex
\hoffset 14truemm
\voffset 31truemm
\font\phead=cmsl9 scaled 950
\font\pnum=cmbx10 scaled 913
\font\pfoot=cmsl9 scaled 950
\headline{\vbox to 0pt{\vskip -4.5mm\line{\small\phead\ifnum
\count0=\startpage ISSN 1472-2739 (on-line) 1472-2747 (printed)
\hfill {\pnum\folio}\else\ifodd\count0\def\\{ }%
\ifx\theshorttitle\relax\thetitle\else\theshorttitle\fi\hfill{\pnum\folio}
\else\def\\{ and }{\pnum\folio}\hfill\ifx\theshortauthors\relax\theauthors
\else\theshortauthors\fi\fi\fi}\vss}}
\footline{\vbox to 0pt{\vglue 0mm\line{\small\pfoot\ifnum\count0=\startpage
\copyright\ \gtp\hfill\else
\agt, Volume \thevolumenumber\ (\thevolumeyear)\hfill\fi}\vss}}
\else
\headsep 23pt
\footskip 35pt
\hoffset -4truemm
\voffset 12.5truemm
\font\lhead=cmsl9 scaled 1050
\font\lnum=cmbx10 
\font\lfoot=cmsl9 scaled 1050
\makeatletter
\def\@oddhead{{\small\lhead\ifnum\count0=\startpage ISSN 1472-2739 
(on-line) 1472-2747 (printed)\hfill {\lnum\number\count0}\else\ifodd\count0
\def\\{ }\ifx\theshorttitle\relax \thetitle \else\theshorttitle\fi\hfill
{\lnum\number\count0}\else\def\\{ and }{\lnum\number\count0}
\hfill\ifx\theshortauthors\relax 
\theauthors\else\theshortauthors\fi\fi\fi}}\def\@evenhead{\@oddhead}
\def\@oddfoot{\small\lfoot\ifnum\count0=\startpage\copyright\ \gtp\hfill\else
\agt, Volume \thevolumenumber\ (\thevolumeyear)\hfill\fi}
\def\@evenfoot{\@oddfoot}
\makeatother
\fi
\let\maketitlepage\makeagttitle
\let\makeshorttitle\maketitlepage
\let\maketitle\maketitlepage


\newwrite\gtoutfile
\long\gdef\makeheadfile{  
{\def\\{, }\def\s{ }
\immediate\openout\gtoutfile head.xxx
\immediate\write\gtoutfile{To: math@arxiv.org}
\immediate\write\gtoutfile{Subject: put OR rep NNNNN:ppppp}
\immediate\write\gtoutfile{--text follows this line--}
\immediate\write\gtoutfile{Proxy-for: \ifx\theasciiauthors\relax
\theauthors\else\theasciiauthors\fi\s<\ifx\theasciiemail\relax\theemail\else\theasciiemail\fi>}
\immediate\write\gtoutfile{\noexpand\\}
\immediate\write\gtoutfile{Authors: \ifx\theasciiauthors\relax
\theauthors\else\theasciiauthors\fi}
{\def\\{ }\immediate\write\gtoutfile{Title: \ifx\theasciititle\relax
\thetitle\else\theasciititle\fi}}
\immediate\write\gtoutfile{Subj-class: GT or SG, GR etc}
\immediate\write\gtoutfile{MSC-class: \theprimaryclass\ifx\thesecondaryclass\relax\else, \thesecondaryclass\fi}
\immediate\write\gtoutfile{Journal-ref: Algebr. Geom. Topol. \thevolumenumber\s
(\thevolumeyear) \startpage-\finishpage}
\immediate\write\gtoutfile{Comments: Published by Algebraic and
Geometric Topology at}
\immediate\write\gtoutfile{\s\s\s  http://www.maths.warwick.ac.uk/agt/AGTVol\thevolumenumber/agt-\thevolumenumber-\thepapernumber.abs.html}
\immediate\write\gtoutfile{\noexpand\\}
\immediate\write\gtoutfile{}
\ifx\theasciiabstract\relax
\immediate\write\gtoutfile{\theabstract}\else
\immediate\write\gtoutfile{\theasciiabstract}\fi
\immediate\write\gtoutfile{}
\immediate\write\gtoutfile{\noexpand\\}
\immediate\write\gtoutfile{}
\immediate\closeout\gtoutfile}}  

\def\maketitlepage{\makeagttitle\makeheadfile}
\let\makeshorttitle\maketitlepage
\let\maketitle\maketitlepage

\lognumber{30}
\volumenumber{1}
\volumeyear{2001}
\papernumber{30}
\pagenumbers{587}{603}
\received{7 April 2001}
\revised{17 October 2001}
\accepted{22 October 2001}
\published{26 October 2001}

\title{Commensurability of graph products}
\authors{Tadeusz Januszkiewicz\\Jacek \'Swi\c atkowski}
\coverauthors{Tadeusz Januszkiewicz\\Jacek \noexpand\'Swi\noexpand\c atkowski}
\asciiauthors{Tadeusz Januszkiewicz\\Jacek Swiatkowski}
\address{Instytut Matematyczny Uniwersytetu 
Wroc\l awskiego\\(TJ: and IM PAN)\\
pl. Grunwaldzki 2/4; 50-384 Wroc\l aw, Poland}
\asciiaddress{Instytut Matematyczny Uniwersytetu 
Wroclawskiego\\(TJ: and IM PAN)\\
pl. Grunwaldzki 2/4; 50-384 Wroclaw, Poland}
\email{tjan@math.uni.wroc.pl, swiatkow@math.uni.wroc.pl}


\abstract
We define graph products of families of pairs of groups
and study the question when two such graph products are commensurable.
As an application we prove linearity of certain graph products.
\endabstract

\primaryclass{20F65}
\secondaryclass{57M07}

\maketitle


\reflist

\key{BH} {\bf M. Bridson, A. Haefliger},
{\it Metric Spaces of Nonpositive Curvature},
Springer, 1999.

\key{B1} {\bf M. Bourdon}, 
{\it Sur les immeubles fuchsiennes et leur type de quasi-isom\'etrie},
Ergod. Th. and Dynam. Sys. 20 (2000),  343-364.

\key{B2} {\bf M. Bourdon}, 
{\it Sur la dimension de Hausdorff de l'ensemble limite d'une familie
de sous-groupes convexes co-compactes},
C. R. Acad. Sci. Paris, t. 325, Serie I (1997), 1097-1100.

\key{D} {\bf M. Davis}, 
{\it Buildings are $CAT(0)$},
in: Geometry and cohomology in group theory (Durham 1994),
Cambridge UP, 1998.

\key{DJ} {\bf M. Davis, T. Januszkiewicz},
{\it Right angled Artin groups are commensurable 
with right-angled Coxeter groups},
J. of Pure and Appl. Algebra 153 (2000), 229-235.

\key{HW} {\bf T. Hsu, D. Wise},
{\it On linear and residual properties of graph products},
Michigan. Math. 46 (1999), 251-259.

\key{Hu}{\bf} {\bf S. P. Humphries},
{\it On representations of Artin groups and the Tits Conjecture}
J. of Algebra 169 (1994), 847-862.

\key{L} {\bf F. T. Leighton},
{\it Finite common coverings of graphs},
J. Comb. Theory (Ser. B) 33 (1982), 231-238.

\key{W} {\bf K. Whyte}, {\it Amenability, Bilipschitz Equivalence,
and the Von Neumann Con\-jecture},  Duke J. Math. 99 (1999), 93-112. 

\endreflist

Graph products are  useful and pretty generalizations of both products and 
free products, intimately linked with right-angled buildings. Part of their 
appeal is their generality: they can be studied in any category with products 
and direct limits.

The question that motivated the present paper was 
``when are the graph products of two families of groups commensurable". 
The inspiration came from a special case considered in [\ref{DJ}] and from a 
conversation with Marc Bourdon on linearity of certain lattices in 
automorphism groups of right-angled buildings.

Here is an answer to the simplest version of this question. 
Recall  first that 
two groups $G, G^\ast$ are {\sl commensurable} if there is a  group $H$
isomorphic to a subgroup of finite index in both $G$ and $G^\ast$;
they are {\sl strongly commensurable} if $H$ has the 
same index in both $G$ and  $G^\ast$.

\proclaim{Theorem 1} 
Let $\Gamma$ be a finite graph, $(G_v)_{v\in V}$, $(G^\ast_v)_{v\in V}$ 
be two families of groups indexed by the vertex set of $\Gamma$.  
Suppose that for every $v\in V$, $G_v$ and $G^\ast_v$ are strongly 
commensurable with the common subgroup $H_v$. 
Then the graph products ${\bf G}=\Pi_{\Gamma}(G_v)_{v\in V}$,
and ${\bf G}^\ast=\Pi_{\Gamma}(G^\ast_v)_{v\in V}$ are strongly 
commensurable: they share a subgroup of index $\Pi_v [G_v:H_v]$.
\endproc

We will prove a slightly more general result
on graph products of pairs of groups.
The proof uses two complementary descriptions
of right-angled building on which a graph product acts.
One of them allows an easy identification of the group acting as 
the graph product, the other allows to compare subgroups.

Theorem 1 and its stronger version formulated in
Section 4 (Corollary 4.2) have several interesting 
special cases discussed in Section 5.

\rk{Acknowledgements}

We would like to thank Marc Bourdon for a conversation 
which inspired this paper, Mike Davis and Jan Dymara for useful comments,
John Meier for directing us to Hsu-Wise paper and \'Swiatos\l aw Gal
for extensive help with the final version of the manuscript.

Both authors were supported by a KBN grant 5 P03A 035 20.

\section{Graph products of pairs}

\rk{Graphs}
A {\it graph} $\Gamma$ on the vertex set $V=V(\Gamma)$ is an antireflexive 
symmetric relation on $V$. Thus our graphs have no loops
and there is at most one undirected edge between two vertices.
Graphs considered in this paper are always finite.

A {\it full subgraph} $\Gamma^+<\Gamma$ on vertices $W\subset V$
is the restriction of the relation to $W$.

A graph is {\it complete} if there is an edge between any two vertices.

A {\it map of graphs} 
$f:\Gamma \to \Gamma^*$ is an injection of sets of vertices
with the property that if there is an edge between $v,w$
then there is an edge between $f(v),f(w)$. 
Thus our maps of graphs are inclusions.

\rk{Graph products}
Let $\Gamma$ be a finite graph, with vertex set $V$.
Suppose for each $v\in V $ one is given a pair of groups $A_v<G_v$. 
For $S$, a complete subgraph of $\Gamma$, define $G_S= 
\Pi_{v\in S}G_v \times \Pi_{v\in V\setminus S}A_v$.
The family of groups $G_S$ together with obvious inclusions
on factors of products gives a direct system of groups
directed by the poset $\cal P$ of complete subgraphs in $\Gamma$,
empty set and singletons included ($G_\emptyset= \Pi_v A_v; 
G_{\{v\}}=G_v\times \Pi_{w\in V\setminus \{v\}}A_w$).

The direct limit of this system
$${\bf G}= \lim\,(G_S)_{S\in{\cal P}} =
\Pi_\Gamma (G_v, A_v)$$
is called {\sl the graph product along $\Gamma$ of the family
of pairs $(G_v, A_v)$.}
To keep notation simple we will denote it for most of the time by
$\bf G$. Note that for $A_v=\{e\}$ we obtain ordinary graph products.

\rk{Graph products are functorial}
If $g: \Gamma\to \Gamma^*$ is a map of graphs, 
and if there is a family of group homomorphisms 
$\omega_v:G_v\to G_{g(v)}^*$, such that $\omega_v(A_v)<A^*_{g(v)}$
then we have induced maps
$\omega_S: G_S\to G^*_{g(S)}$ 
which clearly commute with the maps of direct systems
and consequently induce a homomorphism
$$\omega: {\bf G}\to {\bf G^*}.$$
If $g$ is a surjection on the vertices and 
$\omega_v$ are all surjections, so is the induced homomorphism $\omega$.
If $g$ is an embedding onto a full subgraph and  
$\omega_v$ are injections, so is the induced homomorphism.

\proc{Remark}\key{1.3}%
\rm It follows from functoriality above that if 
$\Gamma$ is a full subgraph of $\Gamma^*$ then graph product of any 
family of pairs along $\Gamma^*$ contains as a subgroup the graph product 
of that family of pairs restricted to $\Gamma$.
In particular, groups $G_S$ inject into $\bf G$.
Thus we can (and will) consider $G_S$ as subgroups of $\bf G$.
\endproc

\rk{Presentations}Graph products can be given in terms of
generators and relations. Suppose that each group $G_v$ is given
by presentation $\langle S_v|R_v  \rangle$ and that $\Sigma_v$
is a set of generators for the subgroup $A_v$ expressed in
terms of generators in $S_v$. Then the graph product
${\bf G}=\Pi_\Gamma(G_v,A_v)$ is given by the presentation
$\langle \cup_{v\in V}S_v| \cup_{v\in V}R_v \cup C  \rangle$,
where $C$ consists of commutators $\{sts^{-1}t^{-1}\}$
whenever $s\in S_v$, $t\in S_w$ and there is an edge between $v$ 
and $w$ in $\Gamma$, or whenever $s\in S_v$, $t\in\Sigma_w$
for some $v\ne w$. 

\rk{Examples}
\items
\item{(1)} Graph product of pairs $\Pi_\Gamma(G_v,A_v)$ along a
complete graph $\Gamma$ is the (direct) product $\Pi_{v\in V}G_v$.
\item{(2)} If $\Gamma$ is an empty graph (i.e. an empty relation
on the vertex set $V$) then the graph product $\Pi_\Gamma(G_v,A_v)$
is the free product of groups 
$G_{\{v\}}=G_v\times \Pi_{w\in V\setminus\{v\}}A_w$
amalgamated along their common subgroup $G_\emptyset=\Pi_{v\in V}A_v$.
\item{(3)} Graph products (with trivial subgroups $A_v$)
of infinite cyclic groups are called
right-angled Artin groups.
\item{(4)} Graph products of cyclic groups of order 2 are called
right-angled Coxeter groups (i.e. Coxeter groups with exponents
2 or $\infty$ only).
\enditems

\eject

\section{The complex $D_{\bf G}$}

\rk{Description of $D_{\bf G}$}
Let $P$ be the realization of the poset $\cal P$ of complete subgraphs 
in $\Gamma$ i.e. the simplicial complex with the vertex set $\cal P$ 
and with simplices corresponding to flags  (i.e. linearly ordered subsets) 
of $\cal P$. For each $S\in{\cal P}$ let
$P_S$ be the subcomplex of $P$ spanned by those vertices $S'\in{\cal P}$
which contain $S$. Note that the poset of subcomplexes $P_S$ with
the reverse inclusion is isomorphic to the poset $\cal P$.
Define a simplicial complex $D_{\bf G}={\bf G}\times P/\sim$ 
where 
the equivalence relation is given by
$(g_1,x_1)\sim(g_2,x_2)$ iff for some $S\in{\cal P}$ we have
$x_1=x_2\in P_S$ and $g_1^{-1}g_2\in G_{S}\subset{\bf G}$.
We denote the point in $D_{\bf G}$ corresponding to a pair $(g,x)\in
{\bf G}\times P$ by $[g,x]$. 
Group $\bf G$ acts on the complex 
$D_{\bf G}$ on the left by $g\cdot[g',x]=[gg',x]$.

One should keep in mind that the complex $D_{\bf G}$ depends on 
the description of the group as a graph product, rather than on 
the group only.

\proclaim{Remark} 
\rm The $\bf G$ action on $D_{\bf G}$ need not be 
effective. Its kernel is the product  $\Pi N_v<\Pi A_v$,
where $N_v$ is the intersection of all $G_v$ conjugates of $A_v$.
Dividing by the kernel of the action is geometrically sound and gives 
the {\sl reduced graph product of pairs}.
For example if all $A_v$ are normal the reduced graph product is 
just the graph product of quotients.
\endproclaim

\rk{Complex of groups ${\bf G}({\cal P})$}
Denote by ${\bf G}(\cal P)$ the simple complex of groups (in the sense
of [\ref{BH}], Chapter II.12) over the poset $\cal P$ defined by the directed
system $(G_{S})_{S\in{\cal P}}$ of groups. In view of the injectivity
discussed in Remark \ref{1.3}, 
Theorems 12.18, 12.20 and Corollary 12.21 of [\ref{BH}] 
imply:

\proc{Proposition}\key{2.4}%
The simplicial complex $D_{\bf G}$ is isomorphic to the development of the 
complex of groups ${\bf G}(\cal P)$ corresponding to the family 
$(i_S)_{S\in{\cal P}}$ of canonical inclusions 
$i_S:G_{S}\to{\bf G}$ into the direct limit. In particular
$D_{\bf G}$ is connected and simply connected.

Moreover the complex of groups associated to the action of $\bf G$
on $D_{\bf G}$ coincides with ${\bf G}(\cal P)$.
\endproc

\rk{$D_{\bf G}$ is a building} 
The complex $D_{\bf G}$ is well known and is sometimes called the 
right-angled building associated to a graph product 
${\bf G}$, see [\ref{D}, Section 5] and [\ref{BH}] (section 12.30 (2)).
It is indeed a Tits building whose appartments are Davis complexes of the 
(right-angled) Coxeter group which is the graph product of $Z_2$'s along 
$\Gamma$. 

\section{Another description of $\bf G$ and $D_{\bf G}$}

\rk{Associated graph product along the complete graph}
Given a finite graph $\Gamma$ on the vertex set $V$ and
a graph product ${\bf G}=\Pi_\Gamma(G_v,A_v)$,
denote by ${\bf G}^c$ the 
graph product of pairs $(G_v,A_v)$ along the complete graph $\Gamma^c$
on the vertex set $V$. 
Put $\omega^c:{\bf G}\to{\bf G}^c$ to be the homomorphism
given by functoriality discussed in Section 1 
and note that $\omega^c$ is surjective.

Let ${\cal P}^c$ be the
poset of complete subgraphs in $\Gamma^c$ (including singletons
and the empty graph) and let $P^c$ be
its realization. The inclusion $\Gamma\to\Gamma^c$ clearly induces
an injective simplicial map $p^c:P\to P^c$ (where $P$ is the realization
of the corresponding poset for $\Gamma$).

\rk{Complex $\Delta_{\bf G}$ and group $\widetilde G$}
Let $D_{{\bf G}^c}$ be the simplicial complex associated to
the graph product ${\bf G}^c$ as in Section 2. 
Denote by $\pi^c:D_{{\bf G}^c}\to P^c$ the simplicial map 
induced by the projection ${\bf G}^c\times P^c\to P^c$.
Put $\Delta_{\bf G}:=
(\pi^c)^{-1}(p^c(P))$ and note that, since the action of 
${\bf G}^c$ on $D_{{\bf G}^c}$ commutes with $\pi^c$,
the subcomplex $\Delta_{\bf G}\subset D_{{\bf G}^c}$ 
is invariant under this action. Thus we will speak about
the (restricted) action of ${\bf G}^c$ on $\Delta_{\bf G}$.
Consider the universal cover $\widetilde{\Delta_{\bf G}}$ 
of $\Delta_{\bf G}$, with
the action of the group $\widetilde G$ which is the  
extension (induced by the covering 
$\widetilde{\Delta_{\bf G}}\to\Delta_{\bf G}$) 
of the group ${\bf G}^c$ by 
the fundamental group $\pi_1(\Delta_{\bf G})$. 

\proc{Theorem}\key{3.3}%
Groups $\widetilde G$ and $\bf G$ are isomorphic,
simplicial complexes $D_{\bf G}$ and 
$\widetilde {\Delta_{\bf G}}$ are equivariantly isomorphic
and the homomorphism $\widetilde G \to {\bf G}^c$ 
induced by the covering $\widetilde{\Delta_{\bf G}}
\to\Delta_{\bf G}$
coincides with the map
$\omega^c:{\bf G} \to {\bf G}^c$. 
\endproc

\prf
Let $f:D_{\bf G}\to \Delta_{\bf G}\subset D_{{\bf G}^c}$ be defined by
$f([g,x])=[\omega^c(g),p^c(x)]$.
This map is easily seen to be surjective and $\omega^c$-equivariant. 
It induces then a morphism $f_*:{\bf G}\backslash\backslash D_{\bf G}\to
{\bf G}^c\backslash\backslash \Delta_{\bf G}$ 
between the complexes of groups
${\bf G}\backslash\backslash D_{\bf G}$ and
${\bf G}^c\backslash\backslash \Delta_{\bf G}$ 
associated to the actions of $\bf G$ on
$D_{\bf G}$ and of ${\bf G}^c$ on $\Delta_{\bf G}$ 
as in [\ref{BH}].

Observe that for a vertex $[g,S]\in D_{\bf G}$ the 
isotropy subgroup of $\bf G$ at $[g,S]$ can be described as 
$\hbox{Stab}({\bf G},[g,S])=gG_Sg^{-1}$.
By substituting $\bf G$ with ${\bf G}^c$ in this observation
we see that the homomorphism 
$\omega^c:{\bf G}\to {\bf G}^c$
maps stabilizers in $D_{\bf G}$ isomorphically to stabilizers
in $D_{{\bf G}^c}$ and hence also in $\Delta_{\bf G}$. The morphism $f_*$
is then isomorphic on local groups. Since moreover the map
between the underlying spaces (quotient spaces of the corresponding
actions) associated to the morphism $f_*$ is a bijection, it follows that
$f_*$ is an isomorphism of complexes of groups.

Let $u:\widetilde{\Delta_{\bf G}}\to \Delta_{\bf G}$ 
be the universal covering map.
As before, by natural equivariance, this map induces a morphism
$u_*:{\widetilde G}\backslash\backslash \widetilde{\Delta_{\bf G}}\to
{\bf G}^c\backslash\backslash \Delta_{\bf G}$ between the complexes of
groups associated to the corresponding actions. 
It follows then from local injectivity of $u$ that the stabilizers
of $\widetilde G$ in $\widetilde{\Delta_{\bf G}}$ are mapped
isomorphically (by the homomorphism $\widetilde G\to{\bf G}^c$ associated
to the covering) to the stabilizers of ${\bf G}^c$ in $\Delta_{\bf G}$,
hence $u_*$ is isomorphic on local groups.
Combining this with equality of the underlying quotient complexes
(which follows directly from the description of $\widetilde G$)
we see that $u_*$ is also an isomorphism of complexes of groups.

Now, since both complexes $D_{\bf G}$ and $\widetilde{\Delta_{\bf G}}$
are connected and simply connected, it follows that they are both
equivariantly isomorphic to the universal covering of the complex
of groups $\Pi_{v\in V}G_v\backslash\backslash \Delta_{\bf G}$ 
acted upon by the fundamental group of this complex of groups. 
Thus the theorem follows. 
\endprf

\rk{Complex $C{\cal X}$}
Consider the family ${\cal X}=(X_v)_{v\in V}$
of quotients $X_v=G_v/A_v$. Denote by $\cal C$ the poset
consisting of all subsets $Y$ in the disjoint union $\cup{\cal X}$
having at most one common element with each of the sets $X_v$. 
We assume that the empty set $\emptyset$ is also in $\cal C$.
Put $C{\cal X}$ to be the realization of the poset $\cal C$
i.e. a simplicial complex with simplices corresponding to linearly
ordered subsets of $\cal C$. 
Alternatively, $C{\cal X}$ is the simplicial cone
over the barycentric subdivision of the join of the family $\cal X$. 

The complex $C{\cal X}$ carries the action of the group $\Pi_{v\in V}G_v$
induced from actions of the groups $G_v$ on the sets $X_v$ (from the left).

\proc{Proposition}\key{3.5}%
The action of ${\bf G}^c$ on the associated 
complex $D_{{\bf G}^c}$ 
is equivariantly isomorphic to the action of $\Pi_{v\in V}G_v$ on $C{\cal X}$.
\endproc

\prf
We will construct a simplicial isomorphism $c:D_{{\bf G}^c}\to C{\cal X}$
as required, defining it first on vertices. Let $[g,S]\in D_{{\bf G}^c}$
be a vertex where $g=\Pi g_v\in\Pi G_v$, $g_v\in G_v$, and $S\subset V$.
Put $$c_0([g,S]):=\{g_vA_v:v\in V\setminus S\}$$ and notice the following
properties:
\items
\item{(1)} for any vertex $[g,S]$ of $D_{{\bf G}^c}$ its image
$c_0([g,S])$ is a well defined vertex in $C{\cal X}$;
\item{(2)} $c_0$ defines a bijection between the vertex sets of
the complexes $D_{{\bf G}^c}$ and $C{\cal X}$;
\item{(3)} both $c_0$ and $c_0^{-1}$ preserve the adjacency relation 
on the vertex sets in the corresponding complexes (where two vertices
are called {\sl adjacent} when they span a 1-simplex).
\enditems

Note that, by definition, both complexes $D_{{\bf G}^c}$ and
$C{\cal X}$ have the following property: each set of pairwise
adjacent vertices in the complex spans a simplex of this complex
(complexes satisfying this property are often called flag complexes).
This property, together with properties (2) and (3) above, imply
that the map $c_0$ induces a simplicial isomorphism 
$c:D_{{\bf G}^c}\to C{\cal X}$.

Now, if $g'=\Pi g_v'\in {\bf G}^c=\Pi G_v$, with $g_v'\in G_v$, we have
$$\eqalign{
g'\cdot c([g,S])&=g'\cdot\{g_vA_v:v\in V\setminus S\}=
\{g_v'g_vA_v:v\in V\setminus S\}\cr
&=c([g'g,S])=c(g'\cdot [g,S]),\cr}
$$
and hence $c$ is equivariant.\endprf

\rk{Alternative description of $\Delta_{\bf G}$}
Denote by $Q$ the quotient of the action of $\Pi_{v\in V}G_v$ on
$C{\cal X}$, and by $q:C{\cal X}\to Q$ the associated quotient map.
$Q$ is easily seen to be the simplicial cone over the barycentric 
subdivision of the simplex spanned by the indexing set $V$ of
the family $\cal X$. Observe now that the equivariant isomorphism
$c:D_{{\bf G}^c}\to C{\cal X}$ of Proposition \ref{3.5} induces an
isomorphism $\varepsilon:P^c\to Q$ of the quotients, and thus we 
have $q\circ c=\varepsilon\circ\pi^c$. In fact $\varepsilon$
is given on vertices by $\varepsilon(S)=V\setminus S$.
Define the map $\delta:P\to Q$ by $\delta:=\varepsilon\circ p^c$. 
Proposition \ref{3.5} implies then the following.

\proc{Corollary}\key{3.6.1}The subcomplex $q^{-1}(\delta(P))\subset 
C{\cal X}$ is invariant under the action of the group $\Pi_{v\in
V}G_v$ and the action of this group restricted to this subcomplex is
equivariantly isomorphic to the action of ${\bf G}^c$ on $\Delta_{\bf
G}$.
\endproc

Slightly departing from the main topic of the paper, we give the
following interesting consequence of Theorem \ref{3.3}.

\proc{Corollary}\key{3.7}%
A graph product (along any finite graph) of pairs $(G_v, A_v)$
is virtually torsion free iff all $G_v$ are virtually torsion free.
\endproc

\prf Since the groups $G_v$ inject into the graph product
${\bf G}=\Pi_\Gamma(G_v,A_v)$, they are clearly virtually torsion free
if their graph product is. To prove the converse, observe that
by Theorem \ref{3.3} $\bf G$ is a semidirect product of the group
${\bf G}^c=\Pi_{v\in V}G_v$ by the fundamental group $\pi_1(\Delta_{\bf G})$.
Since the space $\Delta_{\bf G}$ is finite dimensional and
aspherical (its universal cover $\widetilde{\Delta_{\bf G}}$ is
isomorphic to the Davis' realization of a building, and hence
contractible, see [\ref{D}]), its fundamental group is torsion free
and the corollary follows.
\endprf

\section{Large common subgroups and the proof of Theorem 1}

\rk{Subgroups}
Let $(G_v, A_v)$ and $(G^*_v, A^*_v)$ be two families of pairs 
of groups. Denote by $\bf G$ and $\bf G^*$ the corresponding graph products of 
pairs along the same graph $\Gamma$, and by ${\bf G}^c$ and $({\bf G^*})^c$
the corresponding graph products along the complete graph 
$\Gamma^c$.
Let $\omega^c:{\bf G}\to{\bf G}^c$ and 
$(\omega^*)^c:{\bf G^*}\to({\bf G^*})^c$
be the homomorphisms induced by functoriality from the inclusion map
$\Gamma\to\Gamma^c$.

For each $v\in V$ let $H_v < G_v$ and $H_v^*< G^*_v$ be arbitrary subgroups.
Denote by ${\bf H}$ and $\bf H^*$ preimages
of subgroups $\Pi H_v<\Pi G_v={\bf G}^c$ and 
$\Pi H^*_v<\Pi G^*_v=({\bf G^*})^c$ under
the maps $\omega^c$ and $(\omega^*)^c$ respectively.

\proc{Theorem}\key{4.2}%
If the left actions of $H_v$ on ${G_v/A_v}$
and of $H_v^*$ on ${G^*_v/A^*_v}$ are equivariantly isomorphic
for all $v\in V$ then the actions of $\bf H$ on $D_{\bf G}$ and
of $\bf H^*$ on  $D_{\bf G^*}$ are equivariantly isomorphic.
In particular the subgroups $\bf H$ and $\bf H^*$ are isomorphic.
\endproc

\prf
Let $\cal X$ and ${\cal X}^*$ be the families of the sets of cosets
for the families $(G_v,A_v)$ and $(G_v^*,A_v^*)$ respectively.
Under assumptions of the theorem, the actions of products
$\Pi H_v$ on $C{\cal X}$ and $\Pi H_v^*$ on $C{\cal X}^*$ are
equivariantly isomorphic. Applying Corollary \ref{3.6.1} we conclude that
the actions of the groups $\Pi H_v$ and $\Pi H_v^*$ on the complexes
$\Delta_{\bf G}$ and $\Delta_{{\bf G}^*}$ respectively are
equivariantly isomorphic.

Denote by $\widetilde H$ and $\widetilde H^*$ the preimages
of the products $\Pi H_v$ and $\Pi H_v^*$ by the homomorphisms
$\widetilde G\to\Pi G_v$ and $\widetilde G^*\to\Pi G_v^*$ 
respectively. It follows that the actions of $\widetilde H$ on
$\widetilde\Delta_{\bf G}$ and of $\widetilde H^*$ on 
$\widetilde\Delta_{\bf G^*}$ are equivariantly isomorphic.
But, due to Theorem \ref{3.3}, these actions are equivariantly isomorphic
to the actions of $\bf H$ on $D_{\bf G}$ and of $\bf H^*$ on
$D_{\bf G^*}$ respectively, hence the theorem.
\endprf

\proc{Corollary}\key{4.3}%
Let $(G_v, A_v)$ and $(G^*_v, A^*_v)$ be two 
families of group pairs indexed by the vertex set $V$ of a finite 
graph $\Gamma$.
Suppose that for all $v\in V$ there exist subgroups  $H_v<G_v$ and 
$H^*_v<G^*_v$ of finite index,
such that the left actions of $H_v$ on $G_v/A_v$ and of  $H^*_v$ on 
$G^*_v/A^*_v$ are equivariantly isomorphic. Then the graph products 
${\bf G}=\Pi_\Gamma(G_v,A_v)$ and ${\bf G^*}=\Pi_\Gamma(G_v^*,A_v^*)$ are
commensurable.
\endproc

\prf According to Theorem \ref{4.2} the groups $\bf G$ and $\bf G^*$
share a subgroup ${\bf H}={\bf H^*}$, which is of finite index in
both of them.
\endprf

\proof{Proof of Theorem 1} Under assumptions of Theorem 1
the left actions of the group $H_v$ on $G_v$ and on $G_v^*$ are
clearly equivariantly isomorphic. Then by Corollary \ref{4.3} 
the graph products
$\Pi_\Gamma G_v$ and $\Pi_\Gamma G_v^*$ share a subgroup $\bf H$
which is easily seen to be of index $\Pi_{v\in V}[G_v:H_v]$
in both graph products. 
\endprf

\section{Applications, examples and comments}

\sh{Is strong commensurability a necessary assumption in Theorem 1?} 

Considering free products $Z_2*Z_2$ and $Z_3*Z_3$ shows that one needs a 
hypothesis stronger than commensurability to guarantee commensurability of 
graph products. 
A more delicate example is provided by a family  of graph products along
the pentagon, where at each vertex we put the group $Z_p$.
Bourdon computes in [\ref{B1}] an invariant (conformal dimension at infinity) 
of the hyperbolic groups arising in this way.
His invariant shows that 
as $p$ varies, the graph products are not even quasiisometric,
hence noncommensurable.

A more subtle reason for noncommensurability occurs for free products 
of surface groups. According to Whyte [\ref{W}], the groups
$M_g*M_g$ and $M_h*M_h$ are quasiisometric if $g,h\ge 2$. 
On the other hand, we have the following well known fact.

\proc{Lemma}\key{5.1.1}%
Free products $M_g*M_g$ and $M_h*M_h$ of 
surface groups are not commensurable if $g\ne h$. 
\endproc

\prf
Recall that Kurosh theorem asserts that
if $N$ is a subgroup of finite index $i$ in $L_1*L_2$,
then $N$ is a free product 
$$N_1*N_2*\dots*N_k*F_l,$$ 
where each $N_j$ is a subgroup of finite index in either $L_1$ or $L_2$,
$F_l$ is a free group of rank $l$ and moreover $i=k+l-1$.
Now assume $L_1, L_2$ are fundamental groups of orientable 
aspherical manifolds of the same dimension $m$ (e.g. surface groups). 
One readily sees that $k=b^m(N)=rank H^m(N, Z)$
while $l$ is the rank  of the kernel in $H^1(N, Z)$ 
of the cup product $H^1(N, Z)\times H^{m-1}(N, Z)\to H^m(N, Z)$
interpreted as a bilinear form.
Hence if one knows $N$, one knows the index of $N$
as a subgroup in $L_1*L_2$. This implies that if the free products 
$L_1*L_2$ and $L_1'*L_2'$ of two such group pairs are commensurable
they are strongly commensurable. 

Now, if $g\ne h$ then the groups $M_g*M_g$ and $M_h*M_h$ are not strongly
commensurable, because they have different Euler characteristics.
It follows that these groups are not commensurable.  
\endprf

\sh{Commensurability of graph products as transformation groups}

As it is shown in Section 1, to each graph product $\bf G$ of group pairs
there is associated a right-angled building $D_{\bf G}$ on which $\bf G$
acts canonically by automorphisms. Such buildings corresponding
to different groups $\bf G$ may sometimes be isomorphic. In particular
we have:

\proc{Lemma}\key{5.2.1}%
Let $(G_v,A_v)_{v\in V}$ and $(G_v^*,A_v^*)_{v\in V}$
be two families of groups and subgroups, indexed by a finite set $V$.
Suppose that for each $v\in V$ the indices (not necessarily finite)
$[G_v:A_v]$ and $[G_v^*:A_v^*]$ are equal. Then for any graph $\Gamma$
on the vertex set $V$ the buildings $D_{\bf G}$ and $D_{{\bf G}^*}$
associated to the graph products ${\bf G}=\Pi_\Gamma(G_v,A_v)$
and ${\bf G}^*=\Pi_\Gamma(G_v^*,A_v^*)$ are isomorphic.
\endproc

\prf Observe that, under assumptions of the lemma, the
complexes $D_{{\bf G}^c}$ and $D_{({\bf G}^*)^c}$, and hence also
their subcomplexes $\Delta_{\bf G}$ and $\Delta_{{\bf G}^*}$,
are isomorphic. Since by Theorem \ref{3.3} the buildings $D_{\bf G}$
and $D_{{\bf G}^*}$ are the universal covers of the complexes
$\Delta_{\bf G}$ and $\Delta_{{\bf G}^*}$, the lemma follows.
\endprf

Call two graph products {\sl commensurable as
transformation groups} if their associated buildings are isomorphic
and if they contain subgroups of finite index whose actions
on the corresponding buildings are equivariantly isomorphic.
The arguments we give in this paper show that the graph products 
satisfying our assumptions  are
not only commensurable but also commensurable as transformation groups
(see Theorem \ref{4.2}). 
Closer examination of these arguments shows
that the strong commensurability condition of Theorem 1 (and a more
general condition of Corollary
\ref{4.3}) is not only sufficient, but also necessary for two graph products
of groups (of group pairs respectively) to be commensurable as transformation 
groups. The details of this argument are not completely immediate
but we omit them.

\sh{Special cases of Theorem 1} 
Theorem 1 has interesting special cases resulting from various examples
of strongly commensurable groups. The simplest class of examples is
given by finite groups of equal order. Thus:

\proc{Corollary}\key{5.3.1}%
Let $(G_v)_{v\in V}$ and $(G_v^*)_{v\in V}$ be
two families of finite groups indexed by the vertex set $V$ of a finite
graph $\Gamma$. Suppose that for each $v\in V$ we have
$|G_v|=|G_v^*|$. Then the graph products $\Pi_\Gamma G_v$
and $\Pi_\Gamma G_v^*$ are strongly commensurable.
\endproc

The infinite cyclic group $Z$ and the infinite dihedral group $D_\infty$
are strongly commensurable since they both contain an infinite cyclic
subgroup of index two. Thus a graph product of infinite cyclic groups
(right-angled Artin group)
is commensurable with the corresponding graph product of infinite
dihedral groups which is a right-angled
Coxeter group. 
Thus we reprove a result from [\ref{DJ}]:

\proc{Corollary}\key{5.3.2}%
Right angled Artin groups are commensurable with 
right-angled Coxeter groups.
\endproc

A source of strongly commensurable groups is given by subgroups
of the same finite index in some fixed group. The intersection of
two such subgroups has clearly the same finite index in both of them.
As an example of this kind consider a natural number $g\ge2$ and
a tessellation of the hyperbolic plane $H^2$ by regular $4g$-gons with all 
angles equal to $\pi/2g$ (so that $4g$ tiles meet at each vertex). 
Let $T$ be the group of all symmetries of this tessellation
and $W_g<T$ be the Coxeter group generated by reflections in sides
of a fixed $4g$-gon. Consider also the fundamental group $M_g$
of the closed surface of genus $g$ and note that this group can
be viewed as a subgroup of $T$. Since the groups $W_g$
and $M_g$ have the same fundamental domain in $H^2$ (equal to
a single $4g$-gon) they have clearly the same index in $T$
(equal to $8g$, the number of symmetries of a $4g$-gon) and
hence are strongly commensurable.
Since graph products of Coxeter groups are again Coxeter groups, 
Theorem 1 implies:

\proc{Corollary}\key{5.3.3}Graph products of surface groups are
commensurable with Coxeter groups. \endproc

Pairs of subgroups of the same finite index in a given group
(being thus strongly commensurable) are applied also in the following.

\proc{Proposition}\key{5.3.4}Graph products of arbitrary subgroups 
of finite index in right-angled Coxeter groups are commensurable with
right-angled Coxeter\break groups.
\endproc

\prf 
Since  graph products of right-angled
Coxeter groups remain in this class, 
it is sufficient to show that 
a finite index subgroup
in a right-angled Coxeter group $W$ is strongly commensurable with
another right-angled Coxeter group. This is clearly true for finite
groups, as they are (both groups and subgroups) isomorphic
to products of the group $Z_2$. To prove this for an infinite
group $W$, we will exhibit in $W$ a family $W_n:n\in N$ of subgroups,
indexed by all natural numbers, with $[W:W_n]=n$, 
such that each of the groups $W_n$ is also a right-angled Coxeter group. 

Note that if $W$ is infinite, it contains two generators $t$
and $s$ whose product $ts$ has infinite order in $W$.
Let $D$ be a fundamental domain in the Coxeter-Davis complex
$\Sigma$ of $W$. $D$ is a subcomplex in $\Sigma$ with the 
distinguished set of ``faces", so that reflections with respect to those
faces are the canonical generators of $W$. Since the faces of
the reflections $t$ and $s$ are disjoint, the following complex
$$
D_n:=\cases{D\cup tD\cup stD\cup\dots\cup(st)^kD & if $n=2k+1$ \cr
           D\cup tD\cup stD\cup\dots\cup t(st)^{k-1}D & if $n=2k$ \cr}
$$
is a fundamental domain of a subgroup $W_n<W$ generated by reflections
with respect to ``faces" of this complex. By comparing fundamental
domains we have $[W:W_n]=n$, and the proposition follows.

The algebraic wording of this proof is as follows. An infinite
right angled Coxeter group $(W,S)$ contains an infinite dihedral 
parabolic subgroup$(V,\{s,t\})$. 
The map of $S$ which is the identity on $\{s,t\}$
and sends remaining generators to 1 extends to the homomorphism
$r:W\to V$.The group $V$ contains (Coxeter) subgroups 
generated by $s, (st)^ks(st)^{-k}$ 
and $s, (st)^k t(st)^{-k}$. These have indices $2k, 2k+1$ respectively.
Preimages under $r$ of these subgroups are Coxeter subgroups of $W$
of the same indices.
\endprf

The example discussed just before 
Corollary \ref{5.3.3} generalizes as follows.
Let $(T_v)_{v\in V}$ be a family of topological groups
and let $\Lambda_v\subset T_v$ and $\Lambda_v^*\subset T_v$ 
be two families of lattices such that 
volumes of the quotients  $T_v/\Lambda_v$ and $T_v/\Lambda_v^*$ are
finite and equal for all $v$. Suppose also that for each $v\in V$
there is $t\in T_v$ such that the intersection 
$t^{-1}\Lambda_v t\cap \Lambda_v^*$ has finite index in both 
$\Lambda_v^*$ and the conjugated lattice $t^{-1}\Lambda_v t$.
Then for each $v$ the lattices    $\Lambda_v$ and $\Lambda_v^*$
are strongly commensurable and hence the graph products
$\Pi_\Gamma\Lambda_v$ and $\Pi_\Gamma\Lambda_v^*$ are commensurable
for any graph $\Gamma$ with the vertex set $V$.

For surface groups commensurability condition is a very weak one 
and we have the following: 

\proc{Fact}\key{5.3.5}%
Let $M$ and $N$ be two 2-dimensional orbifolds which are developable. 
Then their fundamental groups
$G_M$ and $G_N$ are strongly commensurable 
iff the orbifold Euler characteristics of $M$ and $N$ are equal.
\endproclaim

Clearly, Fact \ref{5.3.5} allows to formulate the appropriate result
on commensurability of graph products of 2-orbifold groups.
On the other hand,
combining this fact with Theorem 1 and with the argument based on
Kurosh' theorem (as in the proof of Lemma \ref{5.1.1}) one has:

\proc{Corollary}\key{5.3.6}%
Under assumptions and notation of Fact \ref{5.3.5}
the free products $G_M*G_M$ and $G_N*G_N$ are commensurable iff
the orbifold Euler characteristics of $M$ and $N$ are equal.
\endproc

We now pass to applications that require the full strength of Corollary \ref{4.3} 
rather than that of Theorem 1.

\sh{Orthoparabolic subgroups of Coxeter groups}

Recall that parabolic subgroup of a Coxeter group $W$ is the group
generated by a subset $S'$ of the generating set $S$ for $W$.
An {\it orthoparabolic} subgroup of a Coxeter group $W$ is
a normal subgroup $J=\ker\rho$ for a homomorphism $\rho:W\to P$ 
to a parabolic subgroup $P$ such that $\rho|_P=id_P$.
We say that $P$ is the {\it orthogonal parabolic} of $J$.
Note that a homomorphism $\rho$ as above, and hence also
an orthoparabolic subgroup orthogonal to $P$, does not always exist.

Since the left actions of a group $J$ on itself
and on the cosets $W/P$ are equivariantly isomorphic,
Theorem \ref{4.2} implies:

\proc{Corollary}\key{5.4.1}%
If for each $v\in V$ group $J_v$ is an orthoparabolic
subgroup in a Coxeter group $W_v$, orthogonal to a parabolic subgroup $P_v$,
then the graph product $\Pi_\Gamma J_v$ is a subgroup in the graph
product $\Pi_\Gamma(W_v,P_v)$. This subgroup has finite index iff
the subgroups $P_v$ are finite for all $v\in V$.
\endproc

Applying presentations of graph products from Section 1, 
we see that any graph product
$\Pi_\Gamma(W_v,P_v)$ of pairs of a Coxeter group and its parabolic
subgroup is again a Coxeter group. Thus Corollary \ref{5.4.1} implies:

\proc{Corollary}\key{5.4.2}%
A graph product of orthoparabolic subgroups of finite index
in Coxeter groups is a finite index subgroup 
of a Coxeter group.
\endproc

Finite cyclic groups $Z_p$ are orthoparabolic in the dihedral groups 
$D_p$ (as well as $Z$ in $D_\infty$). This again allows to reprove
(and extend) the result of [\ref{DJ}] (compare  \ref{5.3.2} above):

\proc{Corollary}\key{5.4.3}%
Graph products of cyclic groups (among them
right-angled Artin groups) are subgroups of finite index in Coxeter groups.
\endproc

More generally, the {\it even} subgroup of a Coxeter group is the kernel
of the homomorphism $h:W\to Z_2$ which sends all generators of $W$
to the generator of $Z_2$. For example, triangle groups $T(p,q,r)$
and other rotation groups of some euclidean or hyperbolic tessellations
are the even subgroups of the Coxeter reflections groups related
to these tessalations. 
Since these groups are clearly orthoparabolic we have:

\proc{Corollary}\key{5.4.4}%
Graph products of even subgroups of Coxeter
groups are finite index subgroups in Coxeter groups.
\endproc

Although it is fairly hard to find orthoparabolics in general Coxeter 
groups, they are plentiful in right-angled groups, or more generally in 
groups where all entries of the Coxeter matrix are even.
There, for every parabolic subgroup there exist orthogonal to it
orthoparabolics (usually many different ones).

\sh{Graph products of finite group pairs}

Note first that by combining Corollaries \ref{5.4.3} and \ref{5.3.1} we obtain:

\proc{Corollary}\key{5.5.1}%
Graph products $\Pi_\Gamma G_v$ of finite
groups $G_v$ are commensurable with Coxeter groups.
\endproc

Next, applying Corollary \ref{4.3} with trivial groups $H_v$, we have:

\proc{Corollary}\key{5.5.2}%
Graph products $\Pi_\Gamma(G_v,A_v)$
and $\Pi_\Gamma(G_v^*,A_v^*)$ of finite group pairs are commensurable
if  $[G_v:A_v]=[G_v^*:A_v^*]$ for all $v\in V$.
\endproc

An argument referring to above corollaries and using cyclic groups
of orders $[G_v:A_v]$ proves then the following.

\proc{Corollary}\key{5.5.3}%
Graph products of finite group pairs are
commensurable with Coxeter groups.
\endproc

In the rest of this subsection we prove the following slightly
stronger result, under slightly stronger hypotheses:

\proc{Proposition}\key{5.5.4}%
Let $(G_v,A_v)_{v\in V}$ be a family of
pairs of a finite group and its subgroup. Suppose that the left
action of $G_v$ on the cosets $G_v/A_v$ is effective for each $v\in V$.
Then any graph product $\Pi_\Gamma(G_v,A_v)$
is a subgroup of finite index in a Coxeter group.
\endproc

\prf
Canonical action of each of the groups $G_v$ on the cosets $G_v/A_v$
defines a homomorphism $i_v:G_v\to S_{G_v/A_v}=S_{|G_v/A_v|}$ to
the symmetric group on the set of cosets. By the assumption of the
proposition this homomorphism is injective. Consider a subgroup
$\hbox{Stab}(A_v,S_{G_v/A_v})=S_{|G_v/A_v|-1}$ and note that
$i_v(A_v)\subset \hbox{Stab}(A_v,S_{G_v/A_v})$. It follows that
there is a homomorphism $i:\Pi_\Gamma(G_v,A_v)\to 
\Pi_\Gamma(S_{G_v/A_v},\hbox{Stab}(A_v,S_{G_v/A_v}))=
\Pi_\Gamma(S_{|G_v/A_v|},S_{|G_v/A_v|-1})$ between the graph products.
Now for each $v\in V$ the action of $G_v$ on $G_v/A_v$ is easily
verified to be equivariantly isomorphic (by $i_v$) to the
action of the image group $i_v(G_v)$ on the cosets
$S_{G_v/A_v}/\hbox{Stab}(A_v,S_{G_v/A_v})$. It follows from
Theorem \ref{4.2} that the homomorphism $i$ is injective and it
maps the graph product $\Pi_\Gamma(G_v,A_v)$ to the subgroup
of finite index in the graph product 
$\Pi_\Gamma(S_{|G_v/A_v|},S_{|G_v/A_v|-1})$.

Symmetric group $S_{|G_v/A_v|}$ is a Coxeter group and its subgroup 
$S_{|G_v/A_v|-1}$
is a parabolic subgroup. By the remark before Corollary
\ref{5.4.2} a graph product of symmetric group pairs
is a Coxeter group, and thus the proposition follows.
\endprf

\proclaim{Remark}\rm  Removing in Proposition \ref{5.5.4} the assumption
of effectiveness for the actions of $G_v$ on $G_v/A_v$ one can
obtain a similar conclusion for the reduced graph products of
pairs $(G_v,A_v)$ as defined in Section 2. 
\endproclaim

\sh{Groups of automorphisms of locally finite buildings}

It is an open question (except in dimension 1, [\ref{L}]) 
whether any two groups of automorphisms acting
properly discontinuously and cocompactly on a fixed locally finite
right-angled buildings are commensurable as transformation groups.
The building $D_{\bf G}$ associated to a graph product ${\bf G}=
\Pi_\Gamma(G_v,A_v)$
is locally finite iff the indices $[G_v:A_v]$ are finite for all $v\in V$.
The action of $\bf G$ on $D_{\bf G}$ is then
properly discontinuous iff the groups $G_v$ are all finite.
Furthermore, since we always assume that $\Gamma$ is finite, this
action is automatically cocompact.

We may now ask above question in the restricted class of
appropriate graph products.
By using Lemma \ref{5.2.1} and Corollary \ref{5.5.2} we have:

\proc{Corollary}\key{5.6.1}%
Let ${\bf G}=\Pi_\Gamma(G_v,A_v)$ and
${\bf G}^*=\Pi_\Gamma(G_v^*,A_v^*)$ be two graph products of finite group
pairs along the same graph $\Gamma$. Suppose that for each $v\in V$
we have $[G_v:A_v]=[G_v^*:A_v^*]$. Then the associated
buildings $D_{\bf G}$ and $D_{{\bf G}^*}$ are locally finite and isomorphic,
and the actions on them are properly discontinuous and cocompact.
Moreover, the groups $\bf G$ and ${\bf G}^*$ are commensurable
as transformation groups.
\endproc

\proclaim{Remark}\rm  By looking more closely one can show that the
assumptions of Corollary \ref{5.6.1} are necessary for 
the buildings $D_{\bf G}$
and $D_{{\bf G}^*}$ to be locally finite and isomorphic and to
carry properly discontinuous actions of $\bf G$ and ${\bf G}^*$.
Thus the question discussed in this
subsection has positive answer
in the class of (associated actions of) graph products. 
We omit the details of the argument.
\endproclaim

\sh{Linearity of graph products}

In [\ref{DJ}] it was pointed out that commensurability of right-angled 
Artin groups (i.e. graph products of infinite cyclic groups) and 
right-angled 
Coxeter groups implies linearity of the former:
Coxeter groups are linear and groups commensurable with linear groups  
are linear by inducing representation. By the same argument
graph products of groups from various other classes are linear.
For example, Corollaries \ref{5.3.3} and \ref{5.5.3} imply the following.

\proc{Corollary}\key{5.7.1}Graph products of surface groups 
and graph products of pairs of finite groups are linear. 
\endproc

\proclaim{Remark}\rm Bourdon [\ref{B2}] using an entirely different method 
constructed and studied faithful 
linear representations of certain graph products of cyclic groups.
The target of any of his representations is the Lorenz group $SO(N,1)$ and 
the dimension is much smaller than of
ones constructed for that group using Corollary \ref{5.7.1}.
\endproclaim

Without referring to commensurability we still can conclude 
that graph products of any subgroups in Coxeter groups are linear. 
This follows from
the fact that graph products of Coxeter groups are Coxeter groups.
The similar fact for pairs of Coxeter groups and their parabolic
subgroups implies:

\proc{Corollary}\key{5.7.3}Let $(W_v, P_v)$ be a family of pairs
where $W_v$ are Coxeter groups and $P_v$ are their parabolic
subgroups.  For each $v\in V$ let $H_v$ be a subgroup of $W_v$.  Then
any graph product of the family of pairs $(H_v, H_v\cap P_v)$ is a
linear group.
\endproc

\prf A graph product $\Pi_\Gamma(H_v,H_v\cap P_v)$
is a subgroup of $\Pi_\Gamma(W_v,P_v)$ which is a Coxeter group.
\endprf

After this paper was written we've learned from John Meier
about a paper of T. Hsu and D. Wise [\ref{HW}]. There linearity
of graph products of finite groups was established
by embedding them into Coxeter groups. 
Linearity of right-angled Artin groups has been 
proved by S. P. Humphries [\ref{Hu}].

\references

\Addresses
\recd

\bye